\documentclass{amsart}

\begin{document}

\title {On the canonical formula of C. L\'evi-Strauss}
\author{Jack Morava}
\address{Department of Mathematics, Johns Hopkins University, Baltimore,
Maryland 21218}
\email{jack@math.jhu.edu}
\thanks{The author was supported in part by the NSF}
\date {6 June 2003}

\maketitle

\begin{quotation}
`No', was the answer. `We have come to give you metaphors for poetry'. 
W.B. Yeats, {\bf A Vision}
\end{quotation} \bigskip 

\noindent
This note sketches a formal mathematical interpretation for 
the anthropologist Claude L\'evi-Strauss's `canonical formula'
\[
F_x(a):F_y(b) \sim F_x(b):F_{a^{-1}}(y) 
\]
which he has found useful in analysing the structure of myths. Maranda's 
volume [10] is a useful introduction to this subject, which has a somewhat 
controversial history [5] among mathematicians. In view of the 
perspective proposed here (in terms of finite non-commutative groups), 
I believe that skepticism is quite understandable. Nevertheless, I believe
that L\'evi-Strauss knows what he means to say, and that difficulties in 
interpreting his formula are essentially those of translation between the 
languages of disciplines (anthropology and mathematics) that normally 
don't engage in much conversation. \bigskip

\noindent
I am posting this in the hope of encouraging such dialog. This document
is addressed principally to mathematicians, and does not attempt to 
summarize any anthropological background; but in hopes of making it
a little more accessible to people in that field, I have spelled out
a few technical terms in more detail that mathematicians might think
necessary. \bigskip

\noindent
I believe that L\'evi-Strauss perceives the existence of a nontrivial 
anti-automorphism of the quaternion group of order eight; the latter is a 
mathematical object similar to, but more complicated than, the Klein group 
of order four which has appeared elsewhere in his work, cf. eg. [7 \S 6.2 
p. 403], [14 p. 135]. In the next few paragraphs I will define some of 
these terms, and try to explain why I believe the formulation above captures 
what LS means to say. I want to thank the anthropologist Fred Damon for many 
discussions about this topic, and in particular for drawing my attention to 
Maranda's volume, which has been extemely helpful. \bigskip

\noindent
{\bf 1} I was especially happy to find, in at least three separate places 
in that book, anthropologists raising the question of whether the two sides 
of the canonical formula are intended to be understood on a symmetrical 
footing [see [12 p. 35], [13 p. 83], and, most clearly, [3 p. 202]]; in 
other words, whether the symbol $\sim$ in L\'evi-Strauss's formula is 
meant to be an equivalence relation in the mathematical sense. \bigskip

\noindent
If formal terms, an equivalence relation is a relation between two
objects (to be concrete: triangles in the plane), with the properties
that \medskip

\noindent
i) if object $A$ is related to object $B$ (in symbols: $A \sim B$), 
and object $B$ is related to object $C$ (ie $B \sim C$) then necessarily
$A \sim C$, ie object $A$ is related to object $C$; \medskip

\noindent
ii) if $A$ is related to $B$, then $B$ is similarly related to $A$ (this
is the axiom of symmetry, which can be stated symbolically: if $A \sim B$
then $B \sim A$), and finally \medskip

\noindent
iii) the axiom of reflexivity: any object is related to itself, ie $A \sim A$.
\bigskip

\noindent
This notion permits us to disguish equivalence from identity; thus in plane 
geometry the symbol $\sim$ is traditionally used for the relation of 
{\bf similarity}, which means that two triangles have the same angles, but 
are not necessarily of the same size. The point is that things can resemble 
each other in various ways [size, color,\dots], and that there may be good 
reason to compare differing {\bf sorts} of equivalences (in the way that CLS 
says that his formula is about `analogies between analogies'). What C\^ot\'e, 
Racine, and Schwimmer all suggest is an interpretation of the canonical 
formula in which the right-hand side is a {\bf transformation} of the left; 
in more standard mathematical notation, this might be written 
\[
F_x(a):F_y(b) \to F_x(b):F_{a^{-1}}(y) 
\]
The existence of such a transformation turning the left side into the 
right does not preclude that transformation from being an equivalence; all 
it does is allow us to regard axiom ii) above as optional. This fits quite 
naturally with current thinking about category theory. \bigskip

\noindent
{\bf 2} To explain why this is relevant involves a short digression about
the roles of the two {\bf characters} $a,b$ and the two {\bf functions}
$x,y$ in the formula, cf. [13 p. 73-76]. This is in turn related to the role 
of symbolic notation in the first place. \bigskip

\noindent
Mathematicians use what they call {\bf variables}, ie symbols such
as $a,b,x,y$, to express relations which hold for a large class of
objects, cf. [1 Ch. A \S 7]; but to interpret these relations, it is necessary
to understand where the assertion is intended to hold. Thus, for
example, a relation which makes sense for triangles might not make sense
for real numbers. Another important background issue in the interpretation
of mathematical formulas is the role of what are called {\bf quantifiers},
which tell us whether (for example) the formula is intended to hold for
{\bf every} object in an appropriate class, or perhaps only that {\bf some}
object exists, for which the relation holds. \bigskip

\noindent
In the case of the canonical formula, this is particularly important,
because such background information about quantifiers and domain
of validity has been left unspecified. The formula requires one of its 
characters ($a$) to have an associated function ($a^{-1}$), cf. [12 p. 83], 
and it also requires one of the functions ($y$) to play the role of a 
character, on the right side of the formula: this is a key part of the 
formula's assertion, in some ways its central essential double twist. 
The formula is thus intrinsically {\bf un}symmetric: it is not required 
that the character $b$ have an associated function $b^{-1}$, nor that 
the function $x$ have a sensible interpretation as a character. This 
suggests that the canonical formula can be paraphrased as the assertion: 
\bigskip

\noindent
$\bullet$ In a sufficiently large and coherent body of myths we can identify 
characters $a,b$ and functions $x,y$, such that the mythical system 
defines a transformation which sends $a$ to $b$, $y$ to $a^{-1}$, 
and $b$ to $y$, while leaving $x$ invariant. \bigskip

\noindent
This transformation will therefore send the ratio, or formal analogy, $F_x(a):
F_y(b)$ into the ratio $F_x(b):F_{a^{-1}}(y)$: this is the usual statement
of the formula. \bigskip

\noindent
{\bf 3} I don't think this is very controversial; several of the contributors 
to [10] have suggested similar interpretations. I have gone into the question
in some detail, however, to make the point that {\bf if} we can treat the
right-hand side of the canonical formula on an equal footing with the
left-hand side, we should then be able to apply the canonical formula {\bf 
again}; but with $b$ now as the new $a$, $y$ as the new $b$, and $a^{-1}$ as 
the new $y$, defining a chain 
\[
F_x(a):F_y(b) \to F_x(b):F_{a^{-1}}(y) \to F_x(y):F_{b^{-1}}(a^{-1}) \;.
\]
which is consistent with the interpretation of $F_x(a)$ as kind of ratio
$x/a$ of $x$ to $a$: the left and right-hand sides of the chain above 
then become the valid rule
\[
\frac{x/a}{y/b} = \frac{x/y}{b^{-1}/a^{-1}}
\]
for the manipulation of grammar-school fractions. Mosko's [11] variant
\[
F_x(a):F_y(b) \sim F_x(b):F_y(a) \;,
\]
of L\'evi-Strauss's formula also has such an interpretation, when the
algebraic values assigned to the variables lie in a commutative group
in which every element has order two, neutralizing the opposition between
$a$ and $a^{-1}$. This presents Mosko's equation as a version of the CF 
valid in particularly symmetrical situations. \bigskip

\noindent
{\bf 4} I've spelled this argument out because it suggests that interpreting
the canonical formula as expressing the existence of a transformation 
relating its two sides is a useful idea. The remainder of this note will
be concerned with the quaternion group mentioned in the first paragraph,
as an example of a consistent classical mathematical system exemplifying  
L\'evi-Strauss's formula. \bigskip

\noindent
It may be useful to say here that a {\bf group}, in mathematical terminology,
is a system of `elements' (real numbers, for example), together with 
a system of rules for their combination (eg addition). There are lots of
such critters, and some of them are {\bf not} commutative, in the sense
that the order in which we combine the elements may be significant. In 
the case of real numbers, order is not important (and hence it's conceivable
one's checkbook might balance); but rotations in three-dimensional space
(cf. Rubik's cube) form another example of a group, in which the order of 
operations {\bf is} important. \bigskip

\noindent
The quaternion group of order eight (there are other quaternion groups, cf.
[2 \S 5.2]) is the set
\[
Q = \{\pm 1,\; \pm i, \; \pm j, \; \pm k\} \;,
\]
with a noncommutative law of multiplication, in which the product of 
the elements $i$ and $j$ (in that order) is $k$, but the product in
the opposite order is $-k$; in other words,
\[
i \cdot j = k = -j \cdot i, \; j \cdot k = i = -k \cdot j,\; k \cdot i = j 
= -i \cdot k \;.
\]
To complete the `multiplication table' for this group, we have to add the
relations
\[
i \cdot i = j \cdot j = k \cdot k = -1 \;,
\]
as well (last but not least) as the relation $(-1)^2 = +1$. The Klein
group $K$ , which has appeared previously in L\'evi-Strauss's work, can be
similarly described, as a {\bf commutative} version of the group $Q$;
in other words, the multiplication table is as before, except that we
don't bother with the plus and minus signs:
\[
K = \{1,i,j,k\} \;,
\]
given the simpler multiplication table 
\[
i \cdot j = k = j \cdot i, \; j \cdot k = i = k \cdot j, \; k \cdot i = j 
= i \cdot k \;,
\]
together with the relations
\[
i \cdot i = j \cdot j = k \cdot k = 1 
\]
(and of course relations like $1 \cdot i = i = i \cdot 1$, etc.). \bigskip

\noindent
Two groups are {\bf isomorphic} if their elements correspond in a way which 
preserves the multiplication laws: thus in L\'evi-Strauss's writings the
Klein group is described as the set of transformations which send a 
symbol $x$ to the possible values $x,\; -x, \; 1/x, \; -1/x$; the first
such transformation [$x \to x$] corresponds to the `identity' element
$1$ in the presentation of $K$ given above, while the second, ie $x \to -x$,
corresponds to the element $i$; similarly $x \to 1/x$ corresponds to $j$, etc.
It is straightforward to check that the multiplication tables
of these two structures correspond, eg the composition of the transformations
$x \to -x$ (corresponding to $i$) with the composition $x \to 1/x$ 
(corresponding to $j$) is the transformation $x \to 1/(-x) = -1/x$ 
corresponding to $k$, and so forth. \bigskip

\noindent
Similarly, an {\bf anti}-isomorphism of groups is an invertible 
transformation which reverses multiplication: it is a map which sends the 
product of any two elements $g,h$ (in that order) to the product of the 
image elements, in the {\bf reverse} order. In the case of commutative
groups, this is a distinction without a difference, but in the case of
a noncommutative group such as $Q$, it can be significant. \bigskip

\noindent
{\bf 5} For example: the transformation $\lambda: Q \to Q$ which sends
$i$ to $k,\; j$ to $i^{-1} = -i$, and $k$ to $j$ is a nontrivial example
of an antiautomorphism: for example,
\[
\lambda (i \cdot j) = \lambda (k) = j = (-i) \cdot k = \lambda (j) \cdot 
\lambda (i) \;,
\]
while
\[
\lambda (j \cdot k) = \lambda (i) = k = j \cdot (-i) = \lambda (k) \cdot
\lambda (j) \;, 
\]
etcetera. Once this is established, it is easy to check that the assignment
\[
x \mapsto 1, \; a \mapsto i, \; y \mapsto j, \; b \mapsto k
\]
sends the antiautomorphism $\lambda$ to the transformation 
\[
x \to x, \; a \to b, \; y \to a^{-1}, b \to y
\]
defining the canonical formula. \bigskip

\noindent
{\bf 6} Quod, as we say in the trade, erat demonstrandum: this presents
an example of a consistent mathematical system, satisfying a version of
L\'evi-Strauss's formula. \bigskip

\noindent
It is a standard principle of mathematical logic, that the consistency of 
a system of axioms can be verified by giving just {\bf one} example of
an interpretation in which those axioms hold true; but I believe that in 
this case, there may be more to the story. Logicians are concerned with 
questions of logical truth, which can be formulated in terms of the 
commutative group $\{\pm1\}$ [which can alternately be described in 
terms of two-valued `yes-no' judgements, or in terms of the even-odd 
distinction among integers]. The Klein group is an interesting kind of 
`double' of this group, with four elements rather than two, and the 
quaternion group takes this doubling process yet one step further. 
Something similar seems to occur in the study of kinship structures [15], 
but the groups encountered in that field remain necessarily 
commutative. \bigskip

\noindent
{\bf 7} I believe the interpretation proposed here is also helpful in
understanding another aspect of the canonical-formula problem, which 
other commentators have also found confusing: in [8 Ch. 6 p. 156], 
L\'evi-Strauss invokes the formula
\[
F_x(a):F_y(b) \sim F_y(x):F_{a^{-1}}(b) \;.
\]
This differs from the previous version: now $x$ on the left of the
equation becomes $y$ on the right, while $a$ on the left becomes
$x$ on the right, $y$ is transformed into $a^{-1}$, and finally 
$b$ remains invariant. In the framework of paragraph six above,
the assignment
\[
x \mapsto i, \; y \mapsto j, \; a \mapsto k, \; b \mapsto 1 
\]
expresses this transformation as another anti-automorphism of $Q$, 
defined now by 
\[
\sigma(i) = j, \; \sigma (j) = k^{-1} = -k, \; \sigma (k) = i \;.
\]
The two transformations differ by the cyclic transformation
\[
\tau: i \mapsto j \mapsto k \mapsto i
\]
which group-theorists call an outer automorphism, of order three,
of the quaternion group $Q$: in these terms, $\lambda = \tau \circ 
\sigma$. The point is that the symmetries [17] of $Q$ form a larger group;
if we include anti-automorphisms among them, we get a very interesting
group of order twenty-four, in which both transformations $\lambda$
and $\sigma$ might be understood as playing a distinguished role. \bigskip

\noindent
{\bf 8} Perhaps it will be useful to mention that modern mathematical logic 
(cf. eg [6]) is very sophisticated, and is willing to study systems with 
`truth-values' in quite general commutative groups, in a way entirely 
consistent with Chris Gregory's Ramusian [4] precepts; but to my knowledge, 
logic with values in {\bf non}-commutative groups has been studied only 
in contexts motivated by higher mathematics (see eg [9]). This may explain, 
to some extent, the difficulty people have had, in finding a mathematical 
interpretation of L\'evi-Strauss's ideas; but it seems clear to me that 
such an interpretation does exist, and that as far as I can see, it fits 
integrally with L\`evi-Strauss's earlier work on the subject. \bigskip 

\noindent
I hope those who read this will not be offended if I close with a personal 
remark. When I first encountered L\'evi-Strauss's formula, my reaction 
was bemusement and skepticism; I took the question seriously, in large
part because I was concerned that it might represent an aspect of some kind
of anthropological cargo-cult, based on a fetishization of mathematical
formalism. I am an outsider to the field, and can make judgements of 
L\'evi-Strauss's arguments only on the basis of internal consistency (in
so far as I am competent to understand them); but I have to say that
I am now convinced that the man knows his business.

\bibliographystyle{amsplain}

\begin{thebibliography}{99}
\bibitem [13]{13} E. Schwimmer, Is the canonic formula useful in cultural 
description, in [9], p. 56 - 96
\end{thebibliography}

\end{document}